\documentclass[12pt]{amsart}
\usepackage{fullpage,amssymb,amsfonts,amsmath}
\usepackage{enumerate}
\usepackage[shortlabels]{enumitem}
\usepackage{comment}
\usepackage{hyperref} 
\usepackage[capitalise]{cleveref}
\usepackage{url}
\usepackage{todonotes}

\newtheorem{thm}{Theorem}

\title{On the Montgomery-Odlyzko method regarding gaps between zeros of the zeta-function}

\author[Goldston]{Daniel A. Goldston}
\address{San Jose State University}
\email{daniel.goldston@sjsu.edu}
 
\author[Trudgian] {Timothy S. Trudgian}

\address{The University of New South Wales} 
\email{t.trudgian@adfa.edu.au}

\author[Turnage-Butterbaugh]{Caroline L. Turnage-Butterbaugh}
\address{Carleton College}
\email{cturnageb@carleton.edu}

\subjclass[2010]{ 
}

\date{}
\begin{document}

\begin{abstract}
\noindent
Assuming the Riemann Hypothesis, it is known that there are  infinitely many consecutive pairs of zeros of the Riemann zeta-function within 0.515396 times the average spacing. This is obtained using the method of Montgomery and Odlyzko. We prove that this method can never find infinitely many pairs of  consecutive zeros within 0.5042 times the average spacing.
\end{abstract}
\maketitle

\section{Introduction}

The existence of Landau--Siegel zeros (or the Alternative Hypothesis) implies that there are long ranges where all the zeros of the Riemann zeta-function are always spaced no closer than one half of the average spacing. Numerical evidence, however, strongly agrees with the GUE model that suggests there is a positive proportion of consecutive zeros within any small multiple of the average spacing, a conclusion that is also a consequence of Montgomery's pair correlation conjecture. There are three methods in the literature used to study small spacings between zeros of the zeta-function (see \cite{MO} and \cite{CGG}, \cite{Monty}, and \cite{STTB}.) The Montgomery--Odlyzko (M-O) method \cite{MO} produces superior results, albeit under assumption of the Riemann Hypothesis (RH). Nevertheless, we are interested in how far one can push this method.

Let us state the problem more precisely. Write the nontrivial zeros of the Riemann zeta-function $\zeta(s)$ as $\rho = \beta + i\gamma$, where $\beta \in (0,1)$ and $\gamma \in \mathbb{R}$. Let $0 < \gamma_1 \le \gamma_2 \le \cdots \le \gamma_n \le \cdots$ denote the ordinates of the nontrivial zeros of $\zeta(s)$ in the upper half-plane. Since
$$N(T)= \sum_{0 < \gamma \le T}1 \sim \frac{T}{2\pi}\log T,$$ it follows that the gap between consecutive zeros $\gamma_{n+1}-\gamma_n$ is $2\pi/\log{\gamma_n}$ on average. To examine how far gaps deviate from the average, define
\[
\mu = \liminf_{n\to \infty}\frac{\gamma_{n+1}-\gamma_n}{2\pi/\log \gamma_n} \quad \text{and} \quad \lambda = \limsup_{n\to \infty}\frac{\gamma_{n+1}-\gamma_n}{2\pi/\log \gamma_n}.
\]
Trivially, we have that $\mu \le 1 \le \lambda$, and it is expected that $\mu =0$ and $\lambda=\infty$. After much work, the best current result for small gaps under RH is $\mu \le 0.515396$ by Preobrazhenski\u{i} \cite{Pre} and for large gaps under RH is $\lambda\ge3.18$ by  Bui and Milinovich \cite{BM18}. We refer the reader to \cite{CTB} for the history of and progress on this problem.

The result of \cite{Pre} is obtained by an argument based on a method introduced by Montgomery and Odlyzko \cite{MO}. Define, for $T\ge 2$,  $c>0$, and $a_k$ a sequence of complex numbers,
\begin{equation}\label{gully} 
h(c) = c - \frac{\displaystyle \Re \sum_{kn \leq y} a_{k} \overline{a_{kn}} g(n)\frac{\Lambda(n)}{n^{1/2}}}{\displaystyle \sum_{k\leq y} |a_{k}|^{2}},
\end{equation}
where 
\begin{equation*}\label{gdef}
g(n) =  \frac{2\sin \left( \frac{\pi c \log n}{\log T}\right)}{\pi \log n}, 
\end{equation*}
and $y = T^{1-\delta}$ for some small $\delta >0$. Montgomery and Odlyzko proved that if $h(c)>1$ for all sufficiently large $T$ for some choice of $a_k$'s, $c$, and a small $\delta$, then assuming RH we have $\mu \leq c$. For large gaps, if we have $h(c) <1$ then  $\lambda\geq c$. Conrey, Ghosh, and Gonek \cite{CGG} showed that, for any choice of $a_k$
\begin{equation}\label{CGG.5} h(c) <1 \qquad \text{if}\qquad  c<1/2,\end{equation}
which shows that the Montgomery--Odlyzko method is unable to obtain $\mu<1/2$. Due to the connection to Landau--Siegel zeros, it is a tantalizing hope that we might nevertheless reach this barrier. We prove, however, that the Montgomery--Odlyzko method falls well short of being able to prove $\mu\le 1/2$.   Thus a new idea is needed to make further progress on this problem.

\begin{thm}\label{duck}
If $c< 0.5042$, then $h(c)<1$.
\end{thm}

We note that a much weaker version of this result has been known to the experts for some time via unpublished work of the first-named author. We also mention the following information concerning limitations of the Montgomery-Odlyzko method for large gaps between zeros. Conrey, Ghosh, and Gonek \cite[p. 423]{CGG} showed that $h(c) >1$ if $c\ge 6.2$, whence, the Montgomery-Odlyzko method cannot prove the existence of gaps at least 6.2 times the average spacing. In a note added in the proof stage of their paper, Conrey, Ghosh, and Gonek remark that 6.2 may be replaced by 3.74. Correcting for a misprint in their paper, their first result is based on the inequality
\begin{equation}
    h(c) \ge c - 2\left( \frac{c}{\pi}\int_{0}^{1}\frac{|\sin \pi cv|}{v}\, dv\right)^{1/2}.
\end{equation}
Using \emph{Mathematica} one finds that $h(c)>1$ for $c\ge 5.5602\ldots$. Their second improvement result can be obtained from the inequality
\begin{equation}
    h(c) \ge c - 2\left( \frac{c}{\pi}\int_{0}^{\pi c}\left(\frac{\sin v}{v}\right)^2\, dv\right)^{1/2}
\end{equation}
proved by a small change in the proof of the previous bound. One now finds with \emph{Mathematica} that $h(c)>1$ if $c\ge 3.6747\ldots$.

We note that the work by Bui and Milinovich \cite{BM18} uses a different method based on the work of Hall \cite{Hall} and hence is not limited in this way.
 \section{Proof of Theorem \ref{duck}}
We take $0< c <1$. Letting $a_{k} = b_{k}k^{-1/2}$, we obtain from  \eqref{gully} that 
\begin{equation} \label{h} h(c)  \le c + \frac{S}{\sum_{k\le y}\frac{|b_k|^2}{k}}, \quad \text{where} \quad S=   \sum_{kn \leq y} \frac{|b_{k}| |b_{kn}| |g(n)|\Lambda(n)}{kn}. \end{equation}
 
 For any $\alpha$, $\beta >0$ with $4\alpha \beta \ge 1$, we have  $|ab|\le \alpha |a|^2 + \beta |b|^2$, and therefore 
\begin{equation}  \label{1ststep}
|S| \leq \alpha \sum_{kn\leq y} \frac{|b_{k}|^{2}}{k}  |g(n)|\frac{\Lambda(n)}{n} +\beta \sum_{kn\leq y} \frac{|b_{kn}|^{2}}{kn}  |g(n)|\Lambda(n) 
=: \alpha S_1 +\beta S_2. \end{equation}
Using $|\sin x|\le |x|$, we have for $1\le u \le y$ and $0<c<1$ 
\begin{equation} \label{gbound}0< g(u) =  \frac{2\sin \left( \frac{\pi c \log u}{\log T}\right)}{\pi \log u} \le \frac{2c}{\log T} . \end{equation}
To evaluate $S_1$, we
have 
\begin{equation} \label{H(x)}S_1 = \sum_{k\le y}\frac{|b_{k}|^2}{k} H(y/k), \qquad \text{where} \qquad H(x) :=  \sum_{n\le x}g(n)\frac{\Lambda(n)}{n}. 
\end{equation}
Using partial summation with 
\begin{equation}\label{L(x)} L(x) := \sum_{n\le x}\frac{\Lambda(n)}{n} = \log x +O(1),\end{equation}
where the asymptotic formula is elementary, we have 
\begin{equation*}\label{S1eq}
\begin{split}
H(x) &= \int_1^x g(u)dL(u) = L(u)g(u) {\bigg|}_1^x - \int_1^x L(u)g'(u)\, du\\&
=L(x)g(x)-\int_1^x\left(\log u +O(1)\right)g'(u)\, du
\\&  = \Big( g(x)\log x +O(g(x) \Big) - \left( g(x)\log x - \int_1^x\frac{g(u)}{u}\, du + O\left(\int_1^x|g'(u)|\, du\right)\right)
\\& = \int_1^x\frac{g(u)}{u}\, du + O(g(x)) +O\left(\int_1^x|g'(u)|\, du\right).
\end{split} \end{equation*}
By \eqref{gbound} $g(u)\ll 1/\log T$, and since $x\cos x -\sin x \ll x^3$ for $0\le x \ll 1$,
\[ g'(u) =\frac{2}{\pi}\left( 
\frac{\cos \left( \frac{\pi c \log u}{\log T}\right)\frac{\pi c \log u}{\log T} -\sin \left( \frac{\pi c \log u}{\log T}\right)}{ u \log ^2 u}\right)  \ll \frac{\log u}{u(\log T)^3},\]
we have $(\int_1^x|g'(u)|\,du \ll \log^2x/(\log T)^3$, and hence
\[ H(x) = \int_1^x\frac{g(u)}{u}\, du + O\left(\frac{\log^2 x}{(\log T)^3}\right). \]
Thus we conclude, since $y\le T$,
\[ \begin{split}S_1 & =\sum_{k\le y}\frac{|b_{k}|^2}{k}\left(\int_1^{y/k}\frac{g(u)}{u}\, du + O\left(\frac{1}{\log T}\right) \right)\\&
=  \sum_{k\le y}\frac{|b_{k}|^2}{k}  \left(\frac{2}{\pi}\int_{0}^{\frac{\pi c \log (y/k)}{\log T}}  \frac{\sin v}{v}\, dv+ O\left(\frac{1}{\log T}\right) \right) ,\end{split}\]
where we made the change of variable $v= \pi c \log u/\log T$ in the last integral.

For $S_2$ we use \eqref{gbound} and the elementary relation $ \sum_{d|n} \Lambda(d) = \log n$, 
to obtain
\begin{equation*} \begin{split}\label{S2estimate} S_2 &\le \frac{2c}{\log T}\sum_{kn\leq y} \frac{|b_{kn}|^{2}}{kn}\Lambda(n) = \frac{2c}{\log T} \sum_{m\leq y} \frac{|b_{m}|^{2}}{m}\sum_{n|m} \Lambda(n) = \frac{2c}{\log T} \sum_{k\leq y} \frac{|b_{k}|^{2}}{k}\log k .\end{split} \end{equation*}
Hence from \eqref{1ststep} we obtain
\begin{equation} \label{Sbound2} S\le  \sum_{k\le y}\frac{|b_k|^2}{k}\left( \frac{2\beta c \log k}{\log T} + \frac{2\alpha}{\pi} \int_{0}^{\frac{\pi c \log (y/k)}{\log T}}\frac{\sin v}{v}\, dv + O\left(\frac{1}{\log T}\right)\right).
\end{equation}
We
define, for  $1\le w \le y$,  
\begin{equation*} \label{G(x)}  G(w) = G(w,\alpha,\beta,c) := \frac{2\beta c \log w}{\log T} + \frac{2\alpha}{\pi} \int_{0}^{\frac{\pi c \log (y/w)}{\log T}}\frac{\sin v}{v}\, dv ,\end{equation*}
and conclude
\begin{equation} \label{Sbound3} \frac{S}{\sum_{k\le y}\frac{|b_k|^2}{k}} \le  \max_{1\le w\le y} G(w) + O\left(\frac{1}{\log T}\right). \end{equation}
Since $G(w)$ is continuous and differentiable on $[1,y]$, the maximum above exists and occurs at either a critical point or at an endpoint of the interval. It is also clear that the  smallest maximum occurs when $4\alpha \beta = 1$, which we henceforth assume. By the fundamental theorem of calculus
\begin{equation} \label{G'} G'(w) = \frac{2c}{w\log T} \left(\beta - \alpha \ \text{sinc}\left(\frac{\pi c \log (y/w)}{\log T}\right)\right),\quad \text{where}  \quad \text{sinc}(x) := \frac{\sin x}{x}. \end{equation}
\medskip

\noindent{\bf Case 1.} Suppose $\beta\ge \alpha$.  Since $\text{sinc}(x) \le 1$ and $\text{sinc}(x)=1$ if and only if $x=0$, from \eqref{G'} $G(w)$ is increasing on $[1,y]$ and 
$\max G(w) = G(y) = \frac{2\beta c \log y}{\log T} = 2(1-\delta) \beta c \le 2\beta c.$
We chose the smallest value of $\beta$ by taking $ \beta =\alpha =1/2$, which from \eqref{h} and \eqref{Sbound3} recovers \eqref{CGG.5}.\footnote{See the last section for comments on how this approach differs from that of \cite{CGG}.}
\medskip

\noindent{\bf Case 2.} Suppose $\beta < \alpha$. Thus $\beta < 1/2$, and we substitute $\alpha = 1/(4\beta)$.
Since $\text{sinc}\left(\frac{\pi c \log (y/w)}{\log T}\right)$ increases on $w\in[1,y]$, we see $wG'(w)$ decreases through the interval and therefore $G'(w)$ also decreases. Thus there can be at most a single critical point $w=w_0$ where $G'(w_0)=0$. Thus $G(w_0)$ is a relative maximum and the absolute maximum in $[1,y]$.  By \eqref{G'} $w_0$ satisfies
\begin{equation}\label{w0} \text{sinc}\left(\frac{\pi c \log (y/w_0)}{\log T}\right) = \frac\beta\alpha = 4\beta^2, \end{equation}
and, letting $T^{\phi_0} := y/w_0$, this can be written as
\begin{equation} \label{phi0} \text{sinc}(\pi c \phi_0) = 4\beta^2.\end{equation}
Since $G'(y) = 2c(\beta -\alpha)/ (y\log T) <0$ and $G'(1)=  \frac{2c}{\log T} \left(\beta - \alpha \ \text{sinc}(\pi c (1-\delta))\right)$ is positive if $\beta$ is close to $\alpha $ and $c$ is not too small,  it is clear that there are critical points $w_0$, in which case we replace \eqref{Sbound3} with
\begin{equation} \label{Sbound4} \frac{S}{\sum_{k\le y}\frac{|b_k|^2}{k}} \le G(w_0) + O\left(\frac{1}{\log T}\right) \le (1+o(1))G(T^{1-\phi_0 -\delta}). \end{equation}

Using \emph{Mathematica} it is easy to compute the largest $c$ obtainable from \eqref{Sbound4} for which $h(c)<1$. Thus for a given value of $c$ we are seeking the smallest maximum as we vary $\beta$. In performing computations $\delta$ can be taken arbitrarily smaller than the accuracy being used in the calculations, and therefore for computations we can take $\delta=0$ and thus $y=T$ in \eqref{Sbound4}.  We start with an initial choice of $c=c_1=0.5$. Searching with a grid of values of $\beta$ we determine their corresponding values of $\phi_0$ from \eqref{phi0}. We then select a smaller range of $\beta$ containing the smallest maximums of the $G(w_0)$'s found with $h(c_1)<1$, and then replace $c_1$ by a larger value $c_2$ and repeat. This quickly converges. We can stop this process whenever we attain as many digits of accuracy as we desire, at which point we have found values $\beta_n$, $c_n$, $w_n$, and $\phi_n$. We now can check directly that $ h(c_n) \le c_n + G(w,1/4\beta_n,\beta_n,c_n)<1$ for $1\le w\le y$.  In this way we find $c_0=0.5042$, $\beta_0 = 0.476$, $\phi_0= .48025375569\ldots$,  and $h(c_0)\le 0.999993501\ldots$.

\section{A comment on the approach}

In the previous section we recovered the result \eqref{CGG.5} of \cite{CGG} in the simple case that $\beta \ge \alpha$. The proof of \eqref{CGG.5} in \cite{CGG} is different, which we describe here for the interested reader. There the authors use \eqref{gbound} 
and \eqref{L(x)} to obtain

\begin{equation} \label{CGGS1bound} S_1  \le \frac{2c}{\log T} \sum_{k\leq y} \frac{|b_{k}|^{2}}{k}\sum_{n\le y/k} \frac{\Lambda(n)}{n}  = \frac{2c}{\log T} \sum_{k\leq y} \frac{|b_{k}|^{2}}{k} \Big( \log (y/k) +O(1)\Big). \end{equation}
Thus, in place of \eqref{Sbound2} they obtain
\begin{equation*} \label{CGGbound} S \le \frac{2c}{\log T}\sum_{k\leq y} \frac{|b_{k}|^{2}}{k}\big( \alpha \log(y/k) + \beta \log k +O(1)\big).  \end{equation*}
Letting $f(u) = \alpha \log(y/u) + \beta \log u$, one finds that $f(1) = \alpha \log y$,  $f(y) = \beta \log y$, and $f'(u) = \frac{\beta - \alpha}{u}$, and thus  $f(u) \le \max(\alpha, \beta) \log y $ for $1\le u \le y$. The optimal bound is obtained by taking $\alpha=\beta=1/2$, and with this choice
\begin{equation*}\label{CGGresult} S\le  \frac{c \log y +O(1)}{\log T}\sum_{k\leq y} \frac{|b_{k}|^{2}}{k} \le c\sum_{k\leq y} \frac{|b_{k}|^{2}}{k}.\end{equation*}
Substituting into \eqref{h} the authors obtain $h(c)\le 2c $,
and thus \eqref{CGG.5}. Actually in \cite{CGG} the usual choice $\alpha = \beta =1/2$ was used in the argument, which we now see is also the optimal choice when using \eqref{CGGS1bound}.
\section*{Acknowledgements}
Trudgian is supported by ARC DP160100932 and FT160100094; Turnage-Butterbaugh is partially supported by NSF DMS-1901293 and NSF DMS-1854398 FRG. We are grateful to Micah Milinovich for some initial calculations, to Brian Conrey for discussions, to an anonymous reviewer who pointed out a number of improvements in an earlier version of the manuscript, and to Ade Irma Suriajaya for organizing the conference at which we met.

\end{document}